\documentclass{article}
\usepackage[cp1251]{inputenc}
\usepackage[english]{babel}
\usepackage{latexsym,amsfonts,amsthm}

\newtheorem{theorem}{ Theorem}[section]
\newtheorem{lemma}[theorem]{ Lemma}
\newtheorem{definition}[theorem]{ Definition}
\newtheorem{remark}[theorem]{ Remark}

\newtheorem{corollary}[theorem]{ Corollary}
\def\p {\partial }

\input{psfig.sty}
\begin{document}
\large

\begin{center}{\sc Roots of knotted graphs and orbifolds}
\end{center}

  \begin{center} S. Matveev\footnote{Partially supported
by the INTAS  Project "CalcoMet-GT" 03-51-3663  and the RFBR grant
04-01-96014-p204Ural-a}
\end{center}

\vspace{0.5cm}

\section{Introduction}
Let $G$ be a graph (1-dimensional CW complex) in a compact
3-manifold $M$. Following~\cite{Homa}, we will apply to the pair
$(M,G)$ certain simplification moves   as long as possible. What
we get is a {\em root} of $(M,G)$. Our main result is that for any
pair $(M,G)$ the
  root   exists and is unique. Similar results hold for graphs
  with colored edges and for 3-orbifolds, which can be viewed as
  graphs with specific colorings. This generalizes the main result
  of~\cite{petr}. For the  case $G=\emptyset$ (when we are in
  the situation of the Milnor  prime decomposition theorem for 3-manifolds)
  we suggest a   simple  proof  that the   irreducible summands are
   determined uniquely. We begin our exposition  with considering
   this   partial case, since the proof of the main results follows the same lines.

The   paper had been written   during my stay at MPIM Bonn. I
thank  the institute for hospitality and financial support. I thank
C. Petronio  who acquainted me with the problem of spherical splitting of
orbifolds and informed me about a few shortcomings in the first version of this paper.
I thank C.Hog-Angeloni for useful discussions.

\section{Roots of manifolds without graphs}

\begin{definition} \label{df:comressiong} Let   $S$ be a
sphere in the interior of a compact 3-manifold $M$. Then the {\em compression move
} of $M$ along $S$ consists in compressing $S$ to a point and
cutting the resulting
 singular manifold    along that point.
\end{definition}
The image of $S$ under the compression move consists of two
points. Of course, the same result can be obtained by   cutting $M$ along
 $S$ and filling by balls the two  copies of $S$  appearing under the
 cut. If $S$ is trivial (i.e.  bounds a ball), then the
 compression of  $M$ along $S$ produces a copy of $M$ and a
 3-sphere.

\begin{definition} \label{df:rootg} A 3-manifold $R$ is called a {\em root}
of a  3-manifold $M $ if
  $R$ is irreducible  (i.e. contains no nontrivial
spheres) and  can be obtained from $M$ by successive  compressions
along   spheres.
 \end{definition}

 \begin{remark} \label{rmondisg}  Performing   compressions along
  spheres $S_1, S_2, \dots$, we will always assume that each next
 sphere   lies away from the
   two-point  images    of
 the previous spheres.
  Therefore, we may think of  $S_1, S_2, \dots$ as being
 contained in $ M$.
 \end{remark}

\begin{theorem} \label{maing} For any compact 3-manifold   $M$   the
 root   exists and is unique up to homeomorphisms and removing
  connected components homeomorphic to $S^3$.
\end{theorem}
We postpone the proof to Section~\ref{proo}.

\subsection{Kneser's finiteness}
\begin{lemma} \label{boundg} For any compact 3-manifold $M$ there exists
a constant $C_0$ such that
  any sequence of    compression moves
   along
     nontrivial spheres
consists of no more than $C_0$ moves.
\end{lemma}

\begin{proof} We follow the original proof of
H.~Kneser \cite{Kn} (with minor modifications). Let us take $C_0=\beta_1+10t$, where
$\beta_1$ is the dimension over $Z_2$ of $H_1(M;Z_2)$ and $t$ is
the total number  of tetrahedra in a triangulation $T$ of $M$. Let
$n>C_0$.

{\sc Step 1.} Suppose    $S_1, \dots S_n\subset M$ are disjoint
spheres such that all successive compressions along them are
nontrivial. These spheres decompose $M$ into parts called {\em
chambers} such that each $S_i$ corresponds to two boundary spheres
$S^+_i ,S^-_i$ of the chambers. It may happen that both $S^+_i
,S^-_i$ belong to the same chamber. Then we remove $S_i$ from the
sequence  $S_1, \dots, S_n$ and renumber the remaining spheres.
Doing so as long as possible, we get a shorter sequence  $S_1,
\dots, S_m$. Since the total number of removed spheres does not
exceed $\beta_1$, we have $m>10t$. Of course, all successive
compressions along $S_1, \dots S_m$ remain nontrivial. Our profit is that
 now  no
chamber approaches to a sphere   from both sides. It follows
that the following property is true.
\begin{enumerate}
\item[(*)]
 {\em No chamber of the sequence is a punctured ball.}
\end{enumerate}
(Otherwise the compression along the last boundary sphere of a
punctured ball chamber would be trivial).

{\sc Step 2.} We claim that   there exists  another sequence
$S'_1, \dots, S'_m$ consisting of the same number of disjoint
spheres such that the new spheres possess property (*) and are
normal with respect to $T$. The following two observations  are
crucial for the proof.

\begin{enumerate}
\item Any sphere inside a punctured ball decomposes it into two
punctured balls.
\item If a manifold contains a nonseparating sphere, then it is
not a punctured ball.
\end{enumerate}

Let $D$ be a compressing disc for a sphere $S_k, 1\leq k\leq m$,
such that $D\cap (\cap_{i=1}^m S_i)=\partial D\subset S_k$. Denote
by $S_k', S_k''$  two spheres obtained by compressing $S_k$ along
$D$. Let us replace $S_k$ by either $S'_k$ or $S''_k$. It follows
easily from the above  observations that at least one of the
sequences thus obtained satisfies (*). To prove the claim, it
suffices  to recall that any collection of disjoint nontrivial  spheres  can
be normalized by such replacements and  isotopies.

{\sc Step 3.} Let $S'_1, \dots, S'_m$ be   disjoint normal spheres
satisfying (*). They cross each tetrahedron of $T$ along triangle
and quadrilateral pieces called {\em patches}. Let us call a patch
  {\em black}, if it   does not lie between two parallel patches
  of the same type. Each tetrahedron contains at most  $10$ black
  patches: at most 8 triangle patches and at most 2 quadrilateral ones.
   Since $m>10t$, at least one of the spheres  is
  {\em white}, i.e. contains no black patches. Let $C$ be a
  chamber such that $\partial C$ contains a white sphere and a non-white
  sphere. Then $C$ crosses each tetrahedron along some number of
  prisms of the type $P\times I$, where $P$ is a triangle or a
  quadrilateral. Since the patches $P\times \{ 0,1\}$ belong to different spheres,
    $C$ is homeomorphic to $S^2\times
  I$. This contradicts to our assumption that $S'_1, \dots, S'_m$
  satisfy (*).
\end{proof}

\subsection{Proof of Theorem~\ref{maing}}\label{proo}
\begin{definition}\label{complg}   The {\em compression
 complexity} ${\bf c}(M )$ of a compact 3-manifold $M$ is the maximal
possible number of successive  compressions of $M$ along
nontrivial spheres.
\end{definition}

 Lemma~\ref{boundg} shows that ${\bf c}(M )$  is well-defined. Also, it
 follows from the definition   that
    compressions of $M$ along    nontrivial spheres
strictly decrease  ${\bf c}(M )$.

\begin{proof} (Of Theorem~\ref{maing})
To prove the existence,   we compress $ M $ along nontrivial
spheres   as long as possible. Since each  compression strictly
decreases the complexity (which is a nonnegative number), the process is finite
and  the final manifold is a root.

 To
prove the uniqueness,
  assume the converse: suppose that there exists  a compact
  3-manifold  having two different roots.   Among all
 such manifolds  we choose a manifold $ M $ having minimal compression
 complexity.
 Then there exist two sequences of compressions of $ M $ along nontrivial spheres
producing two different roots.  Let the first sequence begin with
 compression along a sphere $S$ while the second   along a sphere $S'$.

{\sc Step 1.} Suppose that $S,S'$ are disjoint.
 Denote by $M_S,M_{S'}$ the manifolds obtained by compressing $M$
 along $S,S'$, respectively. Let $N$ be obtained by compressing
 $M_S$ along $S'$. Of course,
 compression of $ M_{S'} $  along  $S $
 also gives $ N $. Therefore,   $ M_{S} $  and $ M_{S'} $
 have a common root. Indeed,
  one can take any root of $N$.
  On the other hand, inequalities
   $ {\bf c}(M_{S})< {\bf c}(M)$, $ {\bf c}(M_{S'})< {\bf c}(M)$,
  and
   the
  inductive assumption tell us that
 $M_{S}$  and $M_{S'}$ have unique roots. It follows that these roots
  coincide,
 which contradicts to our assumption that they are different.

{\sc Step 2.} Suppose that $S\cap S' $ is nonempty. Using an
innermost circle argument, we compress $S$ along discs contained
in $S'$ as long as possible. This procedure transforms $S$ into a
collection of spheres which intersect neither $S$ nor $S'$. At
least one of those spheres (denote it by $\Sigma$) is nontrivial.
Let us apply  Step 1 twice, to the two   pairs of disjoint nontrivial
spheres $S, \Sigma$ and $\Sigma, S'$.  Clearly, for at least one
case we get a contradiction.
\end{proof}

\section{Roots of knotted graphs}
Now we will consider pairs of the type $(M,G)$, where $M$ is a
compact 3-manifolds and $G$ an arbitrary graph (compact
one-dimensional polyhedron) in $M$.

\subsection{Admissible spheres, compressions, and roots}
\begin{definition} \label{df:admis} A 2-sphere $S$ in $(M,G)$ is called {\em
admissible} if  $S\cap G$ consists of no more than three transverse crossing
points.
\end{definition}
We always assume that an admissible sphere is contained in  the
interior of the manifold.
\begin{definition} \label{df:compressible} An admissible  sphere $S$ in $(M,G)$ is called {\em
compressible} if there is a disc $D\subset M$ such that $D\cap S=\p D$,
 $D\cap G=\emptyset$, and each of the two discs bounded by $\p D$ on $S$
 intersects $G$. Otherwise $S$ is {\em incompressible}.
\end{definition}

\begin{definition} \label{df:trivial} An incompressible sphere $S$ in $(M,G)$ is called {\em
trivial} if it bounds a ball $V\subset M$ such that the pair $(V, V\cap G)$ is
homeomorphic to the pair $(Con (S^2), Con (X))$, where $X\subset S^2$ consists
of $\leq 3$ points and $Con$ is the cone.
  An incompressible nontrivial sphere  is called
{\em essential}.
\end{definition}

\begin{definition} \label{df:comression} Let   $S$ be an incompressible
sphere in $(M,G)$. Then the {\em compression move }
of $(M,G)$ along $S$ consist in
compressing $S$ to a point and cutting the resulting
 singular manifold    along that point.
\end{definition}
 Equivalently, the compression along $S$ can be described as
 cutting $(M,G)$ along
 $S$ and taking disjoint cones over $(S_{\pm}, S_{\pm}\cap G)$, where
 $S_{\pm}$ are two copies of $S$  appearing under the cut.

  If   $(M',G')$ is obtained from $(M,G)$ by compression along
 $S$, we write $(M',G')=(M_S,G_S)$. Note that the image of $S$ under this
 compression consists of two points in $M_S$. We will call them {\em stars}.
 The
 stars lie in $G_S$ if and only if
 $S\cap G\neq \emptyset$.

\begin{figure}
\centerline{\psfig{figure=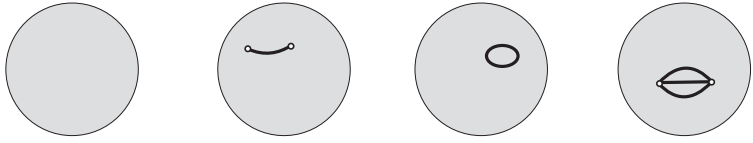,height=1.8cm}}
  \caption{Trivial pairs }
  \label{trivorb}
\end{figure}

\begin{definition} \label{df:trivialpair} A pair $(M,G)$ is called   {\em
trivial} if $M$ is $S^3$ and $G$ is either empty, or a simple arc, or an unknotted
circle, or an unknotted theta-curve  (by an unknotted theta-curve we mean
a graph $\Theta\subset S^3$  such that $\Theta$ is contained in a
disc $D\subset S^3$ and consists of two vertices joined
by three edges).  See Fig.~\ref{trivorb}
\end{definition}

Any trivial pair $(M,G)$ is composed from two copies of a pair $(Con
(S^2),Con(X))$, where $X\subset S^2$ consists of $\leq 3 $ points.

\begin{definition} \label{df:root} A pair $(R,H)$ is called a {\em root}
of a pair $(M,G)$ if:
\begin{enumerate}
\item  $(R,H)$ can be obtained from $(M,G)$ by successive  compressions
along incompressible spheres and
removing   trivial components;
 \item  $(R,H)$
contains no     essential spheres.
 \end{enumerate}
 \end{definition}

 \begin{remark} \label{rmondis}  Performing successive compressions along
 incompressible spheres $S_1, S_2, \dots$, we will always assume that each next
 sphere $S_k$ in the pair $(M_k,G_k)$ obtained from
 $(M_0,G_0)=(M,G)$ by compressing along
 the first $k-1$ spheres lies away from the stars (point  images of
 $S_1,  \dots, S_{k-1}$  under compressions).
  Therefore, we may think of  $S_1, S_2, \dots$ as being
 contained in $(M,G)$.
 \end{remark}

\begin{theorem} \label{ain} For any pair  $(M,G)$   the
 root   exists and is unique up to homeomorphisms and removing
 trivial pairs.
\end{theorem}
We   prove the existence at the end of Section~\ref{existroot} and get
the uniqueness is a corollary of  our main theorem on the uniqueness of
{\em efficient roots}, see Corollary~\ref{noneffroot}.

\subsection{Behavior of spheres with respect to compressions}

We will call a subset $Y$ of $(M,G)$ {\em clean}
if $Y\cap G=\emptyset$. Otherwise $Y$ is {\em dirty}.

\begin{lemma} \label{incomprigid} Let $S,S' $ be disjoint admissible
spheres in $(M,G)$ such that $S$ is incompressible. Then $S'$ is incompressible
in $(M,G)$ if and only if $S'$ is incompressible in $(M_S,G_S)$.
\end{lemma}
\begin{proof} Suppose $S'$ is compressible in $(M,G)$ with compressing disc $D$.
We decrease the number $\#(D\cap S)$ of intersection circles as follows. First,
we choose a disc $D'\subset D$ bounded by an innermost circle
of $D\cap S\subset D$. Since $S$ is incompressible, $\p D'$ bounds a clean disc
$E\subset S$. Then we take a disc $E'\subset E$ bounded by an innermost circle
of $E\cap D\subset E$. Compressing $D$ along $E'$, we get a new compressing
disc $D$ for $S'$ with a smaller number $\#(D\cap S)$. Doing so, we get
a compressing disc $D$ which is  disjoint from $S$ and thus survives
the compression of $(M,G)$ along $S$. It follows that $S'$ remains
compressible. The proof in the other  direction is evident, since we can always
assume that the compressing disc $D\subset (M_S,G_S)$ is away from the stars of
$S$.
\end{proof}

  \begin{lemma} \label{nontrigid} Let $S,S' $ be disjoint incompressible
spheres in $(M,G)$ such that $S'$ is essential in $(M_S,G_S)$.  Then $S'$ is
essential in  $(M,G)$.
\end{lemma}
\begin{proof} By Lemma~\ref{incomprigid}, $S'$ is incompressible
in $(M,G)$. Suppose that  $S'$ is trivial in $(M,G)$.
Then $S$ either does
not intersect the ball $V$ bounded by $S'$ in $(M,G)$  or is a trivial sphere
inside $V$. In both cases $S'$ remains trivial in  $(M_S,G_S)$, a contradiction.
 \end{proof}

 \section{Existence of a root}\label{existroot}

\begin{lemma} \label{bound1} Suppose that $(M,G)$ contains no clean essential
 spheres, i.e. that the manifold $M\setminus G$ is irreducible.
   Then there is a constant $C_1$ depending only on
$(M,G)$ such that any sequence of   compression moves
along essential   spheres consists of no more than $C_1$ moves.
\end{lemma}
\begin{proof} We choose a triangulation $T$ of  $(M,G)$ such
that $G$ is the union of edges and vertices of $T$. Let
$C_1=10t$, where $t$ is the number of tetrahedra in $T$. Consider
a sequence  $S_1,\dots S_n\subset (M,G)$ of $n>C_1$ disjoint spheres such that
each sphere $S_k$ is essential in the pair $(M_k,G_k)$ obtained by compressing
$(M_0,G_0)=(M,G)$ along $S_1,\dots, S_{k-1}$.
 It follows from Lemma~\ref{nontrigid} that
the spheres   are essential in $(M,G)$ and not parallel one to another.

 We claim that there
is a homeomorphism  $h\colon (M,G)\to (M,G)$ such that all spheres $h(S_i), 1\leq
i\leq n$, are normal. Indeed, the usual normalization procedure (see, for
example,~\cite{Ma}) is a
superposition of    moves of two types. The first move is a compression of a sphere along
a   disc inside a triangle face or inside a tetrahedron. The second move
consists in shifting a portion of  $S_i$ along a disc $D\subset M$ such that the
following holds.
\begin{enumerate}
\item The intersection of $D$ with the union of all spheres  is an arc in $\p D\cap
S_i$.
\item The intersection of $D$ with the edges is the complementary arc of $\p D$
contained in the interior of an edge $e$.
\end{enumerate}
Since $M\setminus G$ is irreducible, all moves of the first type
 can be realized by   isotopies
of $(M,G)$. The same is true for the moves of the second type, since
   $e$  cannot  lie in
$G$ (otherwise $S_i$ would be trivial). The terminal homeomorphism
of the normalization isotopy  composed of the above moves is the
required $h$. To prove the lemma, it suffices to apply the same
argument as in Step 3 of the proof of Lemma~\ref{boundg}: since
$n>10t$, there are two spheres $h(S_i), h(S_j)$ such that they bound
$S^2\times I$. This
contradicts to our assumption that all compressions are essential.
\end{proof}

\begin{lemma} \label{bound2} For any pair  $(M,G)$ there exists a constant $C $ such that
  any sequence of    compression moves
   along
  essential  spheres
consists of no more than $C$ moves.
\end{lemma}
\begin{proof} Let $S_1,\dots S_n\subset (M,G)$ be  the given  compression
spheres. We may assume that the pair  obtained by compressions
along all of them   admits no further  compressions along clean
essential spheres. Otherwise we extend the sequence of compression
moves by new compressions along clean essential spheres until
getting a pair with irreducible graph complement.

For any $k, 1\leq k\leq n+1,$ we denote by  $(M_k,G_k)$ the pair
obtained from $(M_0,G_0)=(M,G)$ by compressions along the spheres
$S_1, \dots, S_{k-1}$. Denote also by $(M'_k,G'_k)$ the pair
obtained from $(M_k,G_k)$ by additional compressions along all
remaining clean spheres from the sequence $S_1,\dots S_n$. Then
$(M'_k,G'_k)$ contains no clean essential spheres. It is
convenient to locate the set $X$ of {\em clean stars} (the images
under compressions of all clean essential spheres from $S_1,\dots
S_n$). Then $X$ consists of no more than  $2C_0$ points, where
$C_0=C_0(M,G)$ is the constant from Lemma~\ref{boundg} for a compact 3-manifold
whose interior is $M\setminus G$. We may
think of $X$ as being contained in all $(M'_k,G'_k)$.

Let us decompose the set $S_1,\dots S_n$ into three subsets
$U,V,W$ as follows:
\begin{enumerate}
\item $S_k\in U$ if  $S_k$ is clean.
\item $S_k\in V$ if  $S_k$, considered as a sphere in $(M'_k,G'_k)$,
is an essential  sphere (necessarily dirty).
\item $S_k\in W$ if $S_k$ is a   trivial dirty sphere in $(M'_k,G'_k)$.
\end{enumerate}

Now we estimate the numbers $\# U,\# V,\# W$ of spheres in
$U,V,W$.  Of course, $\# U\leq  C_0$ and $\# V\leq  C_1$, where $C_0$ is
as above and $C_1=C_1(M_0,G'_0)$ is the constant from
Lemma~\ref{bound1}.  Let us prove that $\# W\leq 2C_0$. Indeed,
the compression along each sphere $S_k\subset W$ transforms
$(M'_k,G'_k)$ into a copy of $(M'_k,G'_k)$ and a trivial pair
$(V_k,\Gamma_k)$  containing some number   $w_k$  of clean stars. It is
easy to see that $(V_k,\Gamma_k)$ admits no more than $w_k$
compressions along  essential spheres, and all these spheres are in  $W$.
Since  the total number of clean
stars does not exceed $2C_0$, we get $\# W\leq 2C_0$. Combining
these estimates, we get $n\leq C=3C_0+C_1$.
\end{proof}

\begin{corollary}\label{exists} Any pair $(M,G)$ has a root.
\end{corollary}
\begin{proof}  We apply to $(M,G)$ all possible essential compressions as
long as possible. By Lemma~\ref{bound2} we stop.
\end{proof}

  \section{Efficient roots} \label{eff}
 One of the advantages of   roots introduced above is a flexibility
 of their construction:   each next compression    can be performed
  along any essential sphere.  We pay for that by the  non-uniqueness.
  Indeed,   roots of   $(M,G)$ can   differ by their trivial
  connected components. {\em Efficient roots} introduced in Section~\ref{52} are free from that
  shortcoming.

   \subsection{Efficient systems}

 \begin{definition} A system ${\cal S}=S_1\cup  \dots \cup S_n$
 of disjoint incompressible spheres
  in $(M,G)$ is called {\em efficient} if the following holds:
  \begin{enumerate}
  \item[(1)] compressions along all the spheres give a root of $(M,G)$;
  \item[(2)] any sphere $S_k, 1\leq k\leq n,$ is essential in the
  pair $(M_{{\cal S}\setminus S_k}, G_{{\cal S}\setminus S_k})$ obtained
   from $(M,G)$ by compressions along all
  spheres $S_i,1\leq i\leq n,$ except $S_k$.
\end{enumerate}
  \end{definition}

\begin{figure}
\centerline{\psfig{figure=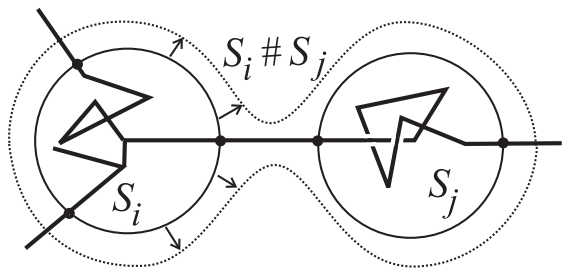,height=3.5cm}}
  \caption{Spherical sliding}
  \label{slide}
\end{figure}

  Evidently, efficient systems exist; to get one, one may construct a system
  satisfying (1) and merely  throw
   away one after another  all spheres not satisfying (2).
Having an efficient system, one can get another one by the following
 moves.
\begin{enumerate}
\item  Let $a\subset (M,G)$ be a clean simple
 arc which joins a sphere $S_i$ with a clean
 sphere $S_j,  i\neq j,$ and has no  common points  with  $\cal S$ except
 its ends. Then the  boundary $\p N$ of a regular neighborhood
 $N(S_i\cup a\cup S_j)$ consists of a copy of $S_i$, a copy of $S_j$, and an
  interior connected sum $S_i\# S_j$ of $S_i$ and $S_j$.  The move consists
 in replacing $S_i$ by $S_i\# S_j$.

 \item The same, but with the following modifications:
 \begin{enumerate}

 \item[i)] $a$ is a simple subarc of $G$ such that all vertices of $G$
 contained in $a$ have valence two, and
 \item[ii)] $S_j$ crosses $G$ in two   points.
\end{enumerate}
\end{enumerate}
 Both moves are called {\em spherical slidings}. See Fig.~\ref{slide}.

\begin{definition} Two efficient system in $(M,G)$ are {\em equivalent} if
one system can be transformed to the other by a sequence of spherical slidings
and an isotopy of $(M,G)$.
\end{definition}

 \begin{lemma}\label{abc} Let  $S$ be an essential sphere in a pair $(U,\Gamma)$  such that
 the pair $(U_S,\Gamma_S)$ obtained by compressing  $(U,\Gamma)$  along
 $S$ is a root, i.e. contains no essential spheres anymore. Suppose  that $D\subset
(U,\Gamma), D\cap S=\p D$ is a compressing disc for $S$ crossing $\Gamma$ in no more than
one transversal point. Additionally we assume that either $D$ is clean or both
discs $D',D''$ into which  $\p D$ decomposes $S$ are dirty.  Denote by $X', X''$ two spheres in $(U,\Gamma)$ obtained by
compressing $S$ along $D$. Then
 there is an isotopy $(U,\Gamma)\to (U,\Gamma)$ taking $S$
 either to   $X'$ or to   $X''$.
\end{lemma}

\begin{proof}   {\em Case 1}. Suppose that $D$ is clean. Since $S$ is
incompressible,  at least one of the spheres $X', X''$ (let $X'$) is
 also clean.
Recall that  $(U_S,\Gamma_S)$ contains no essential spheres. Therefore, $X'$   bounds in $(U_S,\Gamma_S)$
   a clean ball $V$. We denote by $a, a'  $ two stars
  (images of
  $S$ in $(U_S, \Gamma_S)$).

Suppose that $V$, considered as a ball in $ (U_S, \Gamma_S)$, contains
  neither of the two stars $a,a'$.  Then
  we can use $V$   for constructing an isotopy of $S $   to   to $X''$.

   Suppose now that $V$ contains either   $a$ or
   $a'$, but not both. Then the region between $X'$ and $S$ in $(U,\Gamma)$
   is homeomorphic to $S^2\times I$, which assures us that $S$ is isotopic to
   $X'$.

   At last, suppose that $V$ contains both $a,a'$. Then $X''$ is also clean
   and thus bounds a clean ball $W\subset (U_S,\Gamma_S)$. If $W$ contains neither
   $a$ nor  $a'$, or contains only one of them,  then the same arguments
   show that $S$ is isotopic   to
   $X'$ or to $X''$. The case when $W$
    contains both $a,a'$ is impossible, since otherwise $X',X''$
     were parallel and hence
    $S$ were trivial.

 {\em Case 2}.  Suppose that $D$ crosses $\Gamma$ in one point. By assumption, both
 discs   $D', D''$ are dirty. At least one of them (let $D'$)
  crosses  $\Gamma$ in one
 point.
   Since $S$ is incompressible, so is $X'$.
 Then the same argument as in Case 1 shows that $X'$
   bounds in $(U_S,\Gamma_S)$
   a   ball $V$ such that $V\cap \Gamma$ is an unknotted arc.  As above, we
    denote by $a, a'  $ two stars
  (the images of
  $S$ in $(U_S, \Gamma_S)$). Contrary to Case 1, they are points of $\Gamma$ of valence
  two  or three, depending on the number of points in $S\cap \Gamma$.
  Suppose that $V$, considered as a ball in $ (U_S, \Gamma_S)$,
   contains either no stars
    $a,a'$ or   only one of them.    Then
  we can use $V $   for constructing an isotopy of $(U,\Gamma)$ taking $S$ to
  $X'$ or $X''$.
  If  $V$ contains both $a,a'$, then $X''$ is also incompressible,  crosses
  $\Gamma$ in two points and   thus bounds a  ball $W\subset (U_S,\Gamma_S)$
  such that $W\cap \Gamma$ is an unknotted arc. If $W$ contains neither
   $a$ nor  $a'$, or contains only one of them,  then the same arguments
   show that $S$ is isotopic   to
   $X'$ or to $X''$. The case when $W$
    contains both $a,a'$ is impossible, since otherwise $X',X''$
     were parallel and hence
    $S$ were trivial.
   \end{proof}

 \begin{theorem} \label{effic} Any two efficient systems in $(M,G)$ are
 equivalent.
 \end{theorem}

 \begin{proof} Let ${\cal S} $,
 ${\cal S}' $ be two efficient systems in $(M,G)$. Our first
 goal is to replace each system by an equivalent one such that
 the new systems are disjoint.

 {\em Case 1}. Suppose that there is a clean disc $D$ in a sphere
 $S'\subset {\cal S}'$ such that $\p D$ is a circle in a sphere $S\subset {\cal
 S}$ and
 $D\cap {\cal S}=\p D$.
   Let us apply Lemma~\ref{abc} to the pair
 $(U,\Gamma)=
  (M_{{\cal S}\setminus S},G_{{\cal S}\setminus S})$ and  $S,D$ as
  above. We get an isotopy of
  $(U,\Gamma)$ which takes $S$ to
  one of the spheres $X',X''$ (let to $X'$) obtained by
  compressing $S$ along $D$. By construction,
   $\#(X'\cap {\cal S}')<\#(S\cap {\cal S}')$. It is
  easy to see that this isotopy of $S$ in $
  (U,\Gamma)$ can be lifted to
   a composition of isotopies and spherical slidings in $(M,G)$. Each time
   $S$ passes
   a star, we get a spherical sliding.  It means that a new system obtained by
   replacing $S$ by $X'$ is equivalent to   $\cal S$.

 {\em Case 2}. Suppose that all  circles in
 ${\cal S}\cap {\cal S'}$ which are innermost with respect to
 $   {\cal S}$ or to $   {\cal S}'$ bound in  $   {\cal S}$, respectively,
 $   {\cal S}'$ dirty discs. If a sphere from $   {\cal S}$ or  $   {\cal S}'$
 contains at least one circle from ${\cal S}\cap  {\cal S}'$, then it
 contains at least two innermost discs.
 Therefore, at least one of the   discs
 crosses $G$ only once, and we can apply Lemma~\ref{abc} again. As in Case 1,
 this leads us to an equivalent system such that the number of circles in the
 intersection
 is decreased.

 Doing so as long as possible, we get Case 3.

 {\em Case 3}. Suppose ${\cal S},  {\cal S}'$ are disjoint. Our goal is to
 replace $S$ by an equivalent system such that a sphere of $\cal S$ coincides with a
 sphere of ${\cal S}'$.
Since all  spheres of ${\cal S}'$  are trivial in $(M_{\cal S},G_{\cal S})$,
one can
choose an innermost sphere   $S'$.
Then $S'$  bounds a ball $V'\subset (M_{\cal S},G_{\cal S})$
containing a star $a$ of at
least one sphere $S$ of  $  {\cal S}$.
   Note that $V'$   cannot contain   the other star of $S$,
   since otherwise  $S$ would be essential in
  the pair
 $(M_{{\cal S'}}, G_{{\cal S'}})$ obtained from $(M,G)$ by compressions along all
 spheres from $\cal S$. It follows that the spheres $S'=\p V'$ and $S$,
 considered as spheres in $(M_{{\cal S}\setminus S}, G_{{\cal S}\setminus S})$
 are isotopic (see Fig.~\ref{staris}).
 Any isotopy of $S$ to $S'$ can be lifted to $(M,G)$ to
 a composition of
 isotopies and spherical slidings of $S$. The new system $\cal S$ thus obtained
 will have a common sphere $S=S'$ with ${\cal S'}$.

 \begin{figure}
\centerline{\psfig{figure=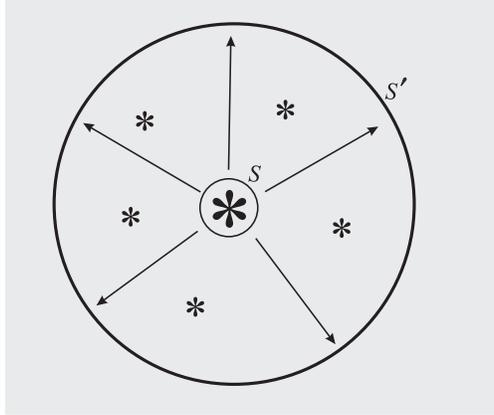,height=5.6cm}}
  \caption{$S$ and $S'$ bound $S^2\times I$ }
  \label{staris}
\end{figure}

  To proceed further, we compress that common sphere and apply the same
  procedure
  to the efficient systems ${\cal S}\setminus S, {\cal S}'\setminus S\subset
  (M_S,G_S)$. We get another pair of coinciding spheres, compress them, apply
  the procedure again, and so on. At the end we get systems consisting of the
  same spheres.\end{proof}

\subsection{Efficient roots} \label{52}
\begin{definition} \label{effroot} A root of $(M,G)$ is {\em efficient}, if
it can be obtained by compressing $(M,G)$ along  spheres of an efficient system.
\end{definition}

 \begin{theorem} \label{main} For any     $(M,G)$   the efficient
 root   exists and is unique up to homeomorphisms.
\end{theorem}
\begin{proof} Evident, since spherical slidings of an efficient
system do not affect the corresponding root.
\end{proof}

\begin{corollary} \label{noneffroot} For any     $(M,G)$   the
 root (not necessarily efficient)     is unique up to homeomorphisms and removing
 trivial pairs.
  \end{corollary}
  \begin{proof} Again evident, since any root can be transformed into an
  efficient one  by removing trivial connected components.
  \end{proof}

 \section{Colored knotted graphs and orbifolds}
 Let $\cal C$ be a set of colors. By a {\em coloring} of a graph $G$ we mean a
 map $\varphi \colon E(G)\to \cal C$, where $E(G)$ is the set of all
 edges of $G$.

 \begin{definition} \label{coloredpair} Let $G_\varphi$ be a colored graph in a
 3-manifold $M$. Then the pair $(M,G_\varphi)$ is called {\em admissible}, if there is
 no incompressible sphere in $(M, G_\varphi)$ which crosses $G_\varphi$ transversely in two points of different
 colors.
 \end{definition}
 It follows from the definition that if $(M,G_\varphi)$ is admissible,
 then $G_\varphi$ has no
 valence two vertices    incident to edges of different colors.
We define   compressions along admissible spheres, trivial pairs,
  roots, efficient systems,   spherical slidings, and efficient roots just in the
 same way as for the uncolored case: we simply forget about the colors.

 \begin{theorem} \label{maincol} For any admissible pair  $(M,G_\varphi )$
   the
 root   exists and is unique up to color preserving homeomorphisms and removing
 trivial pairs. Moreover, any two efficient systems in $(M,G_\varphi )$ are
 equivalent and thus the efficient root is unique up to  color preserving homeomorphisms.
\end{theorem}

\begin{proof}  The proof is literally the same as for the uncolored case. There
are only one place where one should take into account colorings:
  Case 2 of the proof of Lemma~\ref{abc}.
  Indeed,  in this case   there appears
an incompressible  sphere   $X'$  such that
 it crosses $G$  in  two points.
We need to know that
these points have the same colors, and exactly for that purpose one
has imposed the restriction that
the pair
$(M,G_\varphi)$ must be  admissible.     \end{proof}

 Further generalization of the above result consist in
  specifying    sets of {\em allowed}  single colors, pairs of coinciding
  colors, and
  triples of colors. The idea is to allow compressions only along admissible
  spheres whose intersection with $G_\varphi$ belongs to one of the specified
  sets.  Again, all proofs, in particular, the proof of the  corresponding version of
  Theorem~\ref{maincol},  are literally the same with only one exception
   where we need $(M,G_\varphi)$ to be  admissible. We come
  naturally to a generalized version of the {\em orbifold splitting theorem}
  proved recently by C. Petronio~\cite{petr}.

  Recall that a 3-orbifold can be described as a pair $(M,G_\varphi)$, where all
  vertices of $G_\varphi$ have valence 2 or 3 and $G_\varphi$ is colored by the set $\cal C$ of
  all integer numbers
  greater than 1. If $\partial M\neq \emptyset$, then $G\cap \partial M$ should
  consist of univalent vertices of $G_\varphi$. We allow no single colors and allow all pairs of coinciding
  colors.
  The allowed triples are the following: $(2,2,n), n\geq 2$, and $(2,3,k), 3\leq
  k\leq 5$. See~\cite{petr} for background. An orbifold $(M,G_\varphi )$
  is called admissible, if it is admissible in the above sense, i.e.   
if there is
 no incompressible sphere in $(M, G_\varphi)$ which crosses $G_\varphi$ 
 transversely in two points of different
 colors.

  \begin{theorem} \label{orb} For any admissible orbifold  $(M,G_\varphi )$
   the
 root   exists and is unique up to orbifold homeomorphisms and removing
 trivial pairs. Moreover, any two efficient systems in $(M,G_\varphi )$ are
 equivalent and thus the efficient root is unique up to orbifold homeomorphisms.
\end{theorem}

 \end{document}